\documentclass[12pt]{amsart}
\usepackage{amssymb}
\usepackage[margin=1in]{geometry}
\usepackage{hyperref}

\newcommand{\modd}[1]{\,(\textup{mod}\, {#1})}
\def\Re{\textup{Re}\,}

\hypersetup{
  colorlinks= true, 
  urlcolor   = blue, 
  linkcolor  = blue, 
  citecolor  = blue 
}

\theoremstyle{plain}
\newtheorem*{theorem*}{Theorem}
\newtheorem*{lemma*}{Lemma}

\begin{document}

\parskip5pt
\parindent20pt
\baselineskip15pt

\title[Computation of Character Sums]{Harmonic Analysis on the Positive Rationals.  Computation of Character Sums}

\author{P. D. T. A. Elliott}
\address{Department of Mathematics, University of Colorado Boulder, Boulder, Colorado 80309-0395 USA}
\email{pdtae@euclid.colorado.edu}

\author{Jonathan Kish}
\address{Department of Applied Mathematics, University of Colorado Boulder, Boulder, Colorado 80309-0526 USA}
\email{jonathan.kish@colorado.edu}






\maketitle


Let $\mathbb{Q}^*$ denote the multiplicative group of positive rationals and, for fixed integers $a>0$, $A>0$, $b$, $B$ that satisfy $\Delta = aB - Ab \ne 0$, $\Gamma$ its subgroup generated by all but finitely many of the fractions $(an+b)/(An+B)$, $n = 1, 2, \dots$.

In the present paper the authors review their recent publication \cite{elliottkish2017harmonicgroup}, c.f. \href{https://arxiv.org/abs/1602.03263}{arXiv:1602.03263 (2016)}, in which harmonic analysis on $\mathbb{Q}^*$ determines the quotient group $\mathbb{Q}^*/\Gamma$, in particular the positive rationals $r$ that have a product representation
\[
r = \prod_{n_0 < n \le k} \left( \frac{an+b}{An+B} \right)^{\varepsilon_n}, \qquad \varepsilon_n = 0, \pm 1,
\]
the integers $n$ to exceed a given bound, $n_0$, terms in the product with $\varepsilon_n = 0$ understood to have the value 1.

We correct a number of misprints, draw attention to sharpenings implicit in the text and give a detailed account of a corresponding non-intuitive simultaneous, i.e. two dimensional, product representation.  A further representation illustrates consistency  of  the results  with those of the first author's 1985 Grundlehren volume \cite{Elliott1985}.

A major waystation is the following result.\\

\noindent {\bf \cite{elliottkish2017harmonicgroup} Theorem 2}.  \emph{Let integers $a>0$, $A>0$, $b$, $B$ satisfy $\Delta = aB - Ab \ne 0$.  Set $\delta = 6(a,A)(aA)^2\Delta^3$.}

\emph{If a completely multiplicative complex-valued function $g$ satisfies}
\[
g\left( \frac{an+b}{An+B} \right) = c \ne 0
\]
\noindent \emph{on all but finitely many positive integers $n$, then there is a Dirichlet character $\modd{\delta}$ with which $g$ coincides on all primes that do not divide $\delta$.}\\

In the notation of the \emph{Constraints} section of \cite{elliottkish2017harmonicgroup}, define $\alpha = (a,b)$, $\beta = (A,B)$, $a_1 = a\alpha^{-1}$, $b_1 = b\alpha^{-1}$, $A_1 = A\beta^{-1}$, $B_1 = B\beta^{-1}$, $\Delta_1 = a_1B_1 - A_1b_1$, so that $\Delta = \alpha\beta\Delta_1$.  We begin by showing that the value of $\delta$ in \cite{elliottkish2017harmonicgroup} Theorem 2 may be reduced to $6\Delta_1$.\\

\emph{Step 4} of the treatment in \cite{elliottkish2017harmonicgroup} employs the following result.\\

\noindent \emph{\bf \cite{elliottkish2017harmonicgroup} Lemma 6}.  \emph{Let the integers $u_j >0$, $v_j$, $(u_j, v_j) = 1$, $j = 1, 2$, satisfy $\Delta_1 = u_1v_2 - u_2v_1 \ne 0$.  Assume that the primitive Dirichlet character $\chi_D$ satisfies}
\[
\chi_D\left( \frac{ u_1k + v_1}{ u_2k + v_2} \right) = c \ne 0,
\]
\noindent \emph{for all $k$ such that $(u_jk + v_j, D) = 1$, $j = 1, 2$, and that there exists a $k_0$ for which this holds, hence a class $k_0 \modd{D}$.}

\emph{Then $D \mid 6(u_1, u_2)\Delta_1$.}

The argument for this lemma reduces to the case that $D$ is a prime-power, $p^t$, and from the hypothesis derives one of three outcomes:
\begin{enumerate}
\item[$(i)$]  $p^t \mid (u_1, u_2)$;
\item[$(ii)$]  $p^t \mid (u_1v_2 - u_2v_1)$ if $p \ge 3$;
\item[$(iii)$]  $p^{t-1} \mid (u_1v_2 - u_2v_1)$ if $p = 2$ or 3.\\
\end{enumerate}

\emph{The first of these outcomes implies the second, and the factor $(u_1, u_2)$ in the statement of \cite{elliottkish2017harmonicgroup} Lemma 6 is superfluous}.

Note that in the proof of \cite{elliottkish2017harmonicgroup} Lemma 4, the assertion on page 927 that:  \emph{we may set all but one $g(p) = 0$} should read with $g(p) = 1$; or $h(p) = 0$.\\

The subsection \emph{Determination of $G$; practical matters} of \cite{elliottkish2017harmonicgroup}, Section 3, introduces the sums
\[
\eta(\beta, \gamma) = \ell^{-\alpha} \sum_{u \modd{\ell^{\alpha + \max(\beta, \gamma)}}} \chi\left( \frac{a_1u + b_1}{\ell^\beta} \right) \overline{\chi}\left( \frac{A_1u + B_1}{\ell^\gamma} \right)
\]
where $\ell$ is a prime, $\ell^\alpha \| \delta$, $\chi$ is a Dirichlet character $\modd{\ell^\alpha}$, it is understood that $\ell^\beta \mid (a_1u + b_1)$, $\ell^\gamma \mid (A_1u + B_1)$, and $\ell^{\min(\beta, \gamma)} \mid \Delta_1$.

It is noted, \cite{elliottkish2017harmonicgroup} p. 936, that a necessary condition for $\chi$ to be a candidate for the character $g$ on $\mathbb{Q}^*/\Gamma$ that is under consideration is that $|\eta(\beta, \gamma)|$ should have the value of the sum $\eta(\beta, \gamma)$ formed with $\chi$ a principal character.  The following simple result shows that this forces the non-zero summands in $|\eta(\beta, \gamma)|$ to have the same value.

\begin{lemma*} \label{elliottkish2019charactersumslemma}
If $c_j$, $j = 1, \dots, k$ lie in the complex unit disc and satisfy $\left| \sum_{j=1}^k c_j \right| = k$, then the $c_j$ lie on the unit circle and are equal.
\end{lemma*}

\emph{Proof}.  For some real $\theta$, $\sum_{j=1}^k c_j = ke^{i\theta}$.  Then $\sum_{j=1}^k (1 - \Re c_j e^{-i\theta}) = 0$ and the summands are non-negative.

Consideration of the case $\eta(0,0)$ shows that any Dirichlet character $\modd{\ell^\alpha}$ candidate for a character on $\mathbb{Q}^*/\Gamma$ will have a constant value on the reduced fractions $(a_1u + b_1)/(A_1u + B_1)$, $u \modd{l^\alpha}$ and after the updated \cite{elliottkish2017harmonicgroup} Lemma 6, satisfy $\ell^\alpha \mid 6\Delta_1$.

\noindent \emph{Remark}.  A closer reading of \cite{elliottkish2017harmonicgroup} Lemma 6 shows that either $\ell^\alpha \mid A_1$, or $p \le 3$, $p \nmid u_1u_2$ and $p^{t-1} \mid \Delta_1$.

The improvement from $\delta = 6(a,A)(aA)^2\Delta^3$ to $\delta = 6\Delta_1$ in \cite{elliottkish2017harmonicgroup} Theorem 2 follows at once.\\

\noindent \emph{Reduction of simultaneous product representations}

For integers $a>0$, $A>0$, $b$, $B$ that satisfy $\Delta = aB - Ab \ne 0$, let $(\mathbb{Q}^*)^2$ denote the direct product of two copies of $\mathbb{Q}^*$, and $\Gamma$ its subgroup generated by the elements $(an+b)\otimes (An+B)$, $n$ an integer exceeding a given $n_0$.

A character on the quotient group $(\mathbb{Q}^*)^2/\Gamma$ extends to a pair of completely multiplicative functions $g_j$, $j = 1, 2$, on $\mathbb{Q}^*$ that satisfy $g_1(an+b)g_2(An+B) = 1$, $n > n_0$; c.f. Elliott \cite{Elliott1985}, Chapter 15.

The reduction argument given in \cite{elliottkish2017harmonicgroup} Section 5 has a misprint in the last line; a correct version may be found in Elliott \cite{Elliott1985}, Chapter 19.  


We begin with
\[
g_1(an+b) g_2(An+B) = c_1 \ne 0, \quad n > n_0,
\]
where the coefficients $a$, $b$, $A$, $B$ as positive.

Replacing $n$ by $bBn$,
\[
g_1(aBn + 1)g_2(Abn+1) = c_2 \ne 0.
\]

Replacing $n$ by $(aB+1)n+1$
\[
g_1(aBn + 1)g_2(Ab[(aB+1)n+1]+1) = c_3.
\]

Eliminating between these relations,
\[
g_2(Ab(aB+1)n + Ab+1)\overline{g}_2(Abn+1) = c_4 \ne 0.
\]
Note that the corresponding discriminant
\[
\det\left(\begin{array}{cc} Ab(aB+1) & Ab+1 \\ Ab & 1 \end{array} \right)
\]
has the value $Ab\Delta$.  Any Dirichlet character that represents $g_2$ in the manner of \cite{elliottkish2017harmonicgroup} Theorem 2 will be to a modulus $\delta$ that divides $6Ab\Delta$.

The requirement that $b$, $B$ be positive may be obviated as follows.

Replacing the variable $n$ by $n+k$ for a positive integer $k$ moves the coefficient quartet $a, b, A, B$ to $a$, $ak+b$, $A$, $Ak+B$.  This does not affect that values of $(a,b)$, $(A,B)$, or the discriminant.  For all sufficiently large $k$, $\delta \mid 6A(ak+b)\Delta$.  Consecutive values of $k$ show that $\delta \mid 6Aa\Delta$; choosing $k$ to be a multiple of $\delta$, that $\delta \mid 6Ab\Delta$.  Hence $\delta \mid 6A(a,b)\Delta$.

In the notation of \emph{Constraints} $g_2$ is essentially represented by a Dirichlet character to a modulus dividing $6A_1\Delta_1$.

Likewise a Dirichlet character representing $g_1$ will be to a modulus dividing $6a_1\Delta_1$.

\noindent \emph{Remark}.  If $g_j$ is a character on the whole group of rationals and $aA<0$, we may employ $g_j(-1)$ to replace $a$, $A$ by $-a$, $-A$, as necessary.  The outcome of the above argument is then formally the same.

We give an example.\\

\begin{theorem*} \label{elliottkish2019charactersumstheorem}
There are simultaneous representations
\[
a = \prod_{n_0 < n \le k} (5n+1)^{\varepsilon_n}, \quad b = \prod_{n_0 < n \le k} (5n-1)^{\varepsilon_n}, \quad \varepsilon_n = 0, \pm 1
\]
of positive integers $a$, $b$ if and only if $a \equiv 1 \modd{5}$, $(ab)^2 \equiv 1 \modd{5^2}$.

In particular, there are infinitely many simultaneous representations
\[
26 = \prod_{j=1}^k (5n_j + 1)^{\varepsilon_j}, \quad 26 = \prod_{j=1}^k (5n_j - 1)^{\varepsilon_j}
\]
with $\varepsilon_j = \pm 1$, the integers $n_j$ exceeding $n_0$.\\
\end{theorem*}

In this case the general character
\[
g_1 \otimes g_2 : (\mathbb{Q}^*)^2 \to (\mathbb{Q}^*)^2/\Gamma \to \text{unit circle},
\]
coincides on the primes not dividing 30 with a pair of Dirichlet characters to a modulus dividing 300.  We may follow in outline the argument given in the two practical sections of \cite{elliottkish2017harmonicgroup}, the pair $g(an+b)$, $\overline{g}(An+B)$ there replaced by $g_1(an+b)$, $g_2(An+B)$; in effect reduce ourselves to the consideration of
\[
\chi(5n+1)\widetilde{\chi}(5n-1) = c \ne 0
\]
where $\chi$, $\widetilde{\chi}$ are possibly distinct Dirichlet characters to each of the moduli $2^2$, 3, $5^2$.  We consider these cases in turn.\\

\noindent \emph{Modulus 3}.  Choose $t$ to satisfy $5t \equiv 1 \modd{3}$ and set $n = t + 3k$.  Corresponding to $\eta(0,1)$ with $\alpha = 1$, $\beta = 0$, $\gamma = 1$, we see that
\[
\chi(2)\widetilde{\chi}(5k + 3^{-1}(5t-1))
\]
has a constant value when no zero.  Summing over $k \modd{3}$ shows this not to be tenable unless $\widetilde{\chi}$ is principal.

Likewise, $\chi$ is principal.\\

\noindent \emph{Modulus $2^2$}.  Choose $t$ to satisfy $5t \equiv 1 \modd{2^3}$ and set $n = t + 2^3k$.  Corresponding to $\eta(1,3)$ with $\alpha = 2$, $\beta = 1$, $\gamma = 3$, we see that
\[
\chi(2)\widetilde{\chi}(5k + 2^{-3}(5t-1))
\]
has a constant value when not zero.  Summing over $k \modd{2^2}$ shows this to be untenable unless $\widetilde{\chi}$ is principal.

Likewise $\chi$ is principal.

\noindent \emph{Remark}.  The closer reading of \cite{elliottkish2017harmonicgroup} Lemma 6 slightly simplifies these cases.

At this stage, for suitable characters $\chi$, $\widetilde{\chi}$ to the modulus $5^2$, $g_1\overline{\chi}\otimes g_2\overline{\widetilde{\chi}}$ has a constant value on the elements $(5n+1)\otimes(5n-1)$ defining $\Gamma$, and is 1 on all primes save possibly $p = 2, 3$;  by the argument at the end of \cite{elliottkish2017harmonicgroup} Lemma 5, on those also.\\

\noindent \emph{Modulus $5^2$}.  This case is perhaps the most interesting.

If $w$ is a generator of the Dirichlet character group $\modd{5^2}$, there are representations $\chi = w^u$, $\widetilde{\chi} = w^v$ for positive integers $u$, $v$.  The hypothesis that $\chi(5n+1)\widetilde{\chi}(5n-1)$ has a constant value becomes that when not zero
\[
w(( 5n+1)^u (5n-1)^v)
\]
has a constant value.  Since
\[
(5n+1)^u(5n-1)^v \equiv (-1)^v[ 1 + 5(u-v)n] \modd{5^2},
\]
for a suitable non-zero constant $c_0$,
\[
\chi(5n+1)\widetilde{\chi}(5n-1) = c_0w( 5(u-v)n + 1).
\]

The character $w$ is certainly primitive and for a suitable non-zero gaussian sum $\theta$ has a representation
\[
w(s) = \theta^{-1} \sum_{r=1}^{5^2} \overline{w}(r) \exp( 2\pi i 5^{-2} rs).
\]
In particular,
\[
\sum_{n=1}^{5^2} w( 5(u-v)n+1) = \theta^{-1} \sum_{r=1}^{5^2} \overline{w}(r) \exp(2\pi i 5^{-2} r) \sum_{n=1}^{5^2} \exp(2\pi i 5^{-1} r(u-v)n).
\]
The summand $\overline{w}(r) = 0$ unless $5 \nmid r$, in which case the innersum over $n$ is zero unless $5 \mid (u-v)$.

The characters $\chi$, $\widetilde{\chi}$ differ multiplicatively by a character $\modd{5}$.

Without loss of generality we assume $\chi = \widetilde{\chi}\chi_5$ to be such a representation.  Then
\[
1 = g_1(5n+1)g_2(5n-1) = \widetilde{\chi}(5n+1)\widetilde{\chi}(5n-1) = \widetilde{\chi}((5n)^2 - 1);
\]
$\widetilde{\chi}(-1) = 1$.  The character $\widetilde{\chi}$ has order dividing 10; it is the square of an arbitrary character $\modd{5^2}$.

Conversely, such a pair $g_1 = \widetilde{\chi}\chi_5$, $g_2 = \widetilde{\chi}$ satisfies the requirement $g_1(5n+1)g_2(5n-1) = 1$.

This determines the dual group of $(\mathbb{Q}^*)^2/\Gamma$.

The pair $a \otimes b$ belongs to the principal class of $(\mathbb{Q}^*)^2/\Gamma$ if and only if for all pairs $\widetilde{\chi}\chi_5 \otimes \widetilde{\chi}$ with $\widetilde{\chi}$ of order dividing 10, $\widetilde{\chi}(ab)\chi_5(a) = 1$.

Since $\chi_5$ may be any character $\modd{5}$, with $\widetilde{\chi}$ principal $a \equiv 1 \modd{5}$ is required.

Then $\widetilde{\chi}(ab) = 1$ for $\widetilde{\chi}\modd{5^2}$ that satisfy $\widetilde{\chi}(-1) = 1$;  $\chi((ab)^2) = 1$ for all $\chi\modd{5^2}$;  $(ab)^2 \equiv 1 \modd{5^2}$.


In this section we note that any prime $p \ge 5$ dividing the parameter $\delta$ that occurs in Theorem 2 necessarily divides $(a_1, A_1)$, hence $(a,A)$.

If such a prime satisfies $p^r \| \delta$, then $p^r \mid \Delta_1$.  Since any Dirichlet character to a prime-power modulus will induce a character to a modulus that is a higher power of that same prime, without loss of generality we may assume that $p^r \| \Delta_1$.

If now $p \nmid a_1$ and $p^{2r} \mid (a_1n + b_1)$, then $p^r$ divides $A_1(a_1n+b_1)+\Delta_1$, i.e. $a_1(A_1n+B_1)$, hence $A_1n+B_1$.  The terms in the sum $p^r \eta(2r,r)$ are well defined and when non-zero have a constant value.  This remains valid even if the underlying Dirichlet character $\chi \modd{p^r}$ is no longer primitive.

Setting $n = t + p^{2r}k$, where $a_1t + b_1 \equiv 0 \modd{p^{2r}}$, the corresponding summands become
\[
\chi\left( p^{-2r}(a_1t+b_1) + a_1k\right) \overline{\chi}\left( p^{-r} (A_1t+B_1) \right).
\]
Note that $p^{r+1} \nmid (A_1t+B_1)$, otherwise $p^{r+1} \mid \Delta_1$.

Summing over a complete residue class system $k \modd{p^r}$ we see that $\chi$ must be principal.

With emphasis on $\overline{\chi}$ we may draw the same conclusion if $p \nmid A_1$.

A slight modification of the argument shows the conclusion also to be valid if $p = 3$ unless $3 \| \delta$ and $3 \nmid a_1A_1\Delta_1$.

Suppose now that $p^r \| \delta$ and the character $\modd{p^r}$ component of the representing Dirichlet character $\modd{\delta}$ is principal.  Set $\delta_0 = p^{-r}\delta$.  If $\left( (a_1s+b_1)(A_1s+B_1),\delta_0 \right)=1$, then
\[
g\left( \frac{a_1\delta_0 n + a_1s+b_1}{A_1\delta_0n + A_1s+B_1} \right) = \chi_{\delta_0}\left( \frac{a_1s+b_1}{A_1s+B_1} \right)
\]
provided the numerator and denominator of the argument of $g$ are not divisible by $p$.  In particular, c.f. \cite{elliottkish2017harmonicgroup} Lemma 4, second part, $g \overline{\chi}_{\delta_0}(p) = 1$ so long as $p \nmid (a_1, A_1)$.

Without loss of generality we may assume that any prime greater than 3 that does not divide $(a_1, A_1)$ also does not divide $\delta$.

Save for the possible exception noted above, this also holds in the case $p=3$.

In the corresponding simultaneous representation via $(a_1n+b_1) \otimes (A_1n+B_1)$ and employing a pair of Dirichlet characters $\chi_\delta \otimes \chi_\gamma$, $\delta \mid 6a_1\Delta_1$ and has prime factors exceeding 3 only if they divide $a_1$, $\gamma \mid 6A_1\Delta_1$ and has prime factors exceeding 3 only if they divide $A_1$, $\Delta_1 = a_1B_1 - A_1b_1$, as before.

In the example $a_1 = 3$, $b_1 = 1$, $A_1 = 5$, $B_1 = 2$ this permits only that $\chi_\delta$ may be any Dirichlet character $\modd{3}$, $\chi_\gamma$ the principal character $\modd{5}$.  Given $n_0$, there is a simultaneous representation
\[
a = \prod_{n_0 < n \le k} (3n+1)^{\varepsilon_n}, \quad b = \prod_{n_0 < n \le k} (5n+2)^{\varepsilon_n}, \quad \varepsilon = 0, \pm 1,
\]
if and only if $a \equiv 1 \modd{3}$, $(b,5) = 1$.

This accords with the results of \cite{Elliott1985}, Chapter 19.

\bibliographystyle{amsplain}
\bibliography{MathBib}

\end{document}